\documentclass{amsart}
\usepackage{amssymb}
\usepackage{amsmath,amssymb}

\theoremstyle{plain}
\newtheorem{thm}{Theorem}[section]
\newtheorem{theorem}[thm]{Theorem}

\newtheorem{lemma}[thm]{Lemma}

\newtheorem{proposition}[thm]{Proposition}
\theoremstyle{definition}
\newtheorem{remark}[thm]{Remark}

\newtheorem{notation}[thm]{Notation}

\newtheorem{definition}[thm]{Definition}

\newtheorem{example}[thm]{Example}

\numberwithin{equation}{section}

\newcommand{\p}{\partial}

\newcommand{\sA}{{\mathcal A}}

\newcommand{\sC}{{\mathcal C}}
\newcommand{\sD}{{\mathcal D}}
\newcommand{\sE}{{\mathcal E}}
\newcommand{\sF}{{\mathcal F}}

\newcommand{\sH}{{\mathcal H}}
\newcommand{\sI}{{\mathcal I}}

\newcommand{\sL}{{\mathcal L}}

\newcommand{\sO}{{\mathcal O}}

\newcommand{\sR}{{\mathcal R}}

\newcommand{\sT}{{\mathcal T}}

\newcommand{\sV}{{\mathcal V}}
\newcommand{\sW}{{\mathcal W}}

\newcommand{\C}{{\mathbb C}}

\newcommand{\I}{{\mathbb I}}

\newcommand{\BP}{{\mathbb P}}

\newcommand{\U}{{\mathbb U}}
\newcommand{\V}{{\mathbb V}}

\newcommand{\II}{{\rm II}}
\newcommand{\III}{{\rm III}}

\newcommand{\fg}{{\mathfrak g}}

\def\Sym{\mathop{\rm Sym}\nolimits}

\def\Hom{\mathop{\rm Hom}\nolimits}

\title[Lines on contact manifolds]{Lines on holomorphic contact manifolds and a generalization of $(2,3,5)$-distributions to higher dimensions}

\author{Jun-Muk Hwang and Qifeng Li}

\thanks{Jun-Muk Hwang was supported by the Institute for Basic Science (IBS-R032-D1). Qifeng Li was supported by the NSFC grant No. 12201348.}

\begin{document}






\maketitle

\begin{abstract}
Since the celebrated work by Cartan, distributions with \nobreak{small} growth vector $(2,3,5)$  have been studied extensively.
In the holomorphic setting,   there is a natural correspondence  between holomorphic
$(2,3,5)$-distributions  and nondegenerate lines on holomorphic contact manifolds of dimension 5.
We generalize this correspondence to higher dimensions by studying nondegenerate lines on holomorphic contact manifolds and  the corresponding class of distributions of small growth vector $(2m, 3m, 3m+2)$ for any positive integer $m$.
\end{abstract}


\section{Introduction}\label{s.1}
We work in the holomorphic setting. All manifolds and maps are holomorphic, unless stated otherwise. Open subsets refer to Euclidean topology. A Zariski-open subset of a complex manifold is the complement of a closed analytic subset.
We use the following terminology on distributions.

\begin{definition}\label{d.distribution}
A {\em distribution} on a complex manifold $M$ is a vector subbundle $D \subset TM$ of the tangent bundle of $M$. \begin{itemize} \item[(i)] Lie brackets of local sections of $D$ give a homomorphism $${\rm Levi}^D: \wedge^2 D \to TM/D,$$  called the {\em Levi tensor} of $D$.
\item[(ii)] There is a sequence of vector bundles, called the {\em weak derived system} of $D$,  $$D = \p^{(0)}D \subset \p D = \p^{(1)}D \subset \p^{(2)} D \subset \cdots \subset \p^{(d)} D = \p^{(d+1)} D $$ defined on a Zariski-open subset of $M$ such that their associated sheaves satisfy $$\sO(\p^{(i+1)}D) = [ \sO(\p^{(i)}D), \sO(D)] + \sO(\p^{(i)} D)$$ for each $0 \leq i \leq d+1.$ We say that $D$ is {\em regular} at $x \in M$, if the weak derived system of $D$ is a sequence of vector bundles in a neighborhood of $x$.
    \item[(iii)]
  The sequence of integers $$({\rm rank}(D), {\rm rank}(\p D), {\rm rank}(\p^{(2)}D), \cdots, {\rm rank}(\p^{(d)}D))$$ is called the {\em small growth vector}  of $D$.
  \item[(iv)] When $D$ is regular at $x \in M$,
  the graded vector space $${\rm symb}_x(D) := \oplus_{i=1}^d  (\p^{(i)}D)_x / (\p^{(i-1)}D)_x$$
has a natural structure of a nilpotent graded Lie algebra induced by Lie brackets of local sections.
It is called the {\em symbol algebra} of $D$ at $x$. \end{itemize} \end{definition}

\bigskip
In the celebrated paper \cite{Cr}, Cartan studied   distributions  on 5-dimensional manifolds with  the small growth vector $(2,3,5)$,  commonly called $(2,3,5)$-distributions.  Cartan  investigated the local equivalence problem for $(2,3,5)$-distributions and  many mathematicians have developed this study further (see the references in \cite{IKTY}).  A remarkable development in 1990's was the theory of  abnormal extremals   of $(2,3,5)$-distributions, originated from  geometric control theory (see \cite{BH} and \cite{Ze} and the references therein). It  associates a certain contact manifold of dimension 5 to each $(2,3,5)$-distribution.   By reinterpreting this result from the viewpoint of complex geometry, the following one-to-one correspondence has been discovered in Theorems 5.10 and 5.12 of \cite{HL}, where `small growth vector' is abbreviated to s.g.v. :
\begin{equation*} \left\{ \begin{array}{cc} \mbox{distributions of} \\ \mbox{s.g.v.} (2,3,5) \end{array} \right\}  \stackrel{(\textsf{Corr.$1$})}{\Longleftrightarrow}
\left\{ \begin{array}{cc} \mbox{nondegnerate lines on contact} \\  \mbox{manifolds of dimension 5} \end{array} \right\}
 \end{equation*}
See Definition \ref{d.Gauss} for the precise meaning of the right hand side. Note that nondegenerate lines are called `contact unbendable rational curves of Cartan type'  in Definition 5.8 of \cite{HL}. A local version  of this correspondence, where the right hand side is replaced by  Lagrangian cone structures on contact manifolds of dimension 5 satisfying certain conditions, is given in Theorem 3.1 of \cite{IKTY}. From the viewpoint of \cite{HL}, the Lagrangian cone structure is a local description of the VMRT of the nondegenerate lines (see Lemma \ref{l.Legendrianvmrt}).

In the current paper, we generalize this correspondence to higher dimensions as the following one-to-one correspondence:
\begin{equation*} \left\{ \begin{array}{cc} \mbox{some distributions of} \\  \mbox{s.g.v.} (2m, 3m, 3m+2)   \end{array} \right\}   \stackrel{(\textsf{Corr.$m$})}{\Longleftrightarrow}
\left\{ \begin{array}{cc} \mbox{nondegnerate lines on contact} \\  \mbox{manifolds of dimension $2m+3$} \end{array} \right\}
 \end{equation*}
  The major difference from the case $m=1$ is that the distributions on the left hand side are not determined by the small growth vector alone:  their symbol algebras must be of the form $\fg_+(F)$ described in Definition \ref{d.fg+}, where $F$ is a nondegenerate cubic form on a vector space of dimension $m$. For simplicity, we call distributions on the left hand side {\em $\fg_+$-distributions}. Our main results are Theorems \ref{t.symbol} and \ref{t.lines}, which give a precise statement of the correspondence (\textsf{Corr.$m$}).
 When $m=1$, the Lie algebra $\fg_+(F)$ is determined by the small growth vector because all nondegenerate cubic forms on a vector space of dimension 1 are isomorphic. Thus $\fg_+$-distributions in dimension $5$ are just $(2,3,5)$-distributions. There is more than one  isomorphism type of nondegenerate cubic forms when $m \geq 2.$
 So when $m \geq 2$, the small growth vector alone cannot determine the type of our distributions. As a matter of fact, when $m \geq 3$, there are nontrivial moduli of nondegenerate cubic forms and   the isomorphism types of symbol algebras of a $\fg_+$-distribution  may vary from point to point.
A classical example of (\textsf{Corr.$m$}) is the following.

\begin{example}\label{e.EKP}
 Let $\fg$ be a simple Lie algebra and let $G$ be a complex Lie group with Lie algebra $\fg$.
   Let $X^{\fg}$ be the \emph{adjoint variety} of $\fg$, namely, the highest weight orbit of the coadjoint representation on $\BP \fg^{\vee}$, and let $2m+3$ be the dimension of $X^{\fg}$. The variety $X^{\fg}$ has  a natural $G$-invariant contact structure  and is covered by nondegenerate lines.  The space of  lines on $X^{\fg} \subset \BP \fg^{\vee}$ is a rational homogeneous space $Y^{\fg}= G/P$ for some parabolic subgroup $P \subset G$. If $\fg$ is not of type A or C, there exists a  $G$-invariant distribution $D^{\fg} \subset TY^{\fg}$ of rank $2m$ whose symbol algebras are
   isomorphic to $\fg_+(F)$ for  a homaloidal EKP cubic $F$ (using the terminology in Theorem 3 of \cite{Do}),  which is the cubic form whose associated cubic hypersurface is either
\begin{itemize} \item[(i)] the union of a hyperplane and a quadratic cone with one isolated singular point outside the hyperplane, or \item[(ii)] the secant variety of one of the four Severi varieties. \end{itemize}
The correspondence  (\textsf{Corr.$m$}) in this case is :
\begin{equation*}  \begin{array}{cc} \mbox{distribution $D^{\fg}$ on $Y^{\fg}$ of} \\  \mbox{s.g.v.} (2m, 3m, 3m+2)  \end{array}   \Longleftrightarrow
 \begin{array}{cc} \mbox{ lines on the contact manifold } \\   X^{\fg} \mbox{ of dimension } 2m+3.\end{array}
 \end{equation*}

\end{example}

One interesting consequence of (\textsf{Corr.$m$}) is that there are many examples of holomorphic contact manifolds covered by nondegenerate lines, which is already nontrivial when $m=1$. Another interesting aspect of (\textsf{Corr.$m$}) is its potential application in Riemannian geometry.   A well-known  conjecture in complex geometry, which is equivalent to  the LeBrun-Salamon conjecture (p. 110 of \cite{LS}) on quaternionic-K\"ahler manifolds, is that a Fano contact manifold  whose automorphism group is reductive  is biholomorphic to an adjoint variety in Example \ref{e.EKP}.
There is an  approach to this conjecture  using  lines on Fano contact manifolds of Picard number 1 (for example, see \cite{BKK}, \cite{Ke}, \cite{Ke2}).
 The   lines in this case  are expected to be nondegenerate (see Remark \ref{r.LM})  and we have the associated $\fg_+$-distributions via (\textsf{Corr.$m$}).  We believe that  it is important to understand the geometry of these distributions for this approach to the LeBrun-Salamon conjecture. Note that, if the cubic forms  are homaloidal EKP cubics, it is already proved that the Fano contact manifold is an adjoint variety (by Main Theorem in Section 2 of \cite{Mk}).

     Let us discuss briefly  the content of the paper and the methods employed. The key ingredients in establishing (\textsf{Corr.$m$}) come from deformation theory of rational curves, in particular, the theory of VMRT, some standard results of which  are recalled in Section \ref{s.vmrt}.   Section \ref{s.contact} discusses how to go from the right hand side to the left hand side of  (\textsf{Corr.$m$}). The proof of the main result,  Theorem \ref{t.symbol}, is much more involved than the proof in the case of $m=1$, because the structure of the symbol algebras of the distributions is more intricate  in higher dimensions.   Its proof   uses some  special features of deformation theory of rational curves on   contact manifolds.   Section \ref{s.zelenko}  discusses how to go from the left hand side to the right hand side of  (\textsf{Corr.$m$}). The proof of the main result,  Theorem \ref{t.lines},   is a generalization of that of Theorem 5.10 in \cite{HL} to higher dimensions and the  idea has originated from  \cite{Ze}.

            \section{Unbendable rational curves and VMRT}\label{s.vmrt}

\begin{definition}\label{d.urc}
            Let $X$ be a complex manifold of dimension $n$ and let ${\rm Douady}(X)$ be its Douady space parameterizing all compact analytic subspaces of $X$.
            \begin{itemize}
\item[(i)]    A smooth rational curve $\BP^1 \cong C \subset X$ is {\em unbendable} if its normal bundle $N_C$ is isomorphic to $\sO(1)^{\oplus p} \oplus \sO^{\oplus(n-1-p)}$ for some nonnegative integer $p \leq n-1$. Consequently,
    $$TX|_C \cong \sO(2) \oplus \sO(1)^{\oplus p} \oplus \sO^{\oplus(n-1-p)}$$
    with $TC \subset TX|_C$ corresponding to $\sO(2)$. In this case, we denote by $N_C^+$ (resp. $TX|_C^+$) the subbundle of $N_C$ (resp. $TX|_C$) corresponding to $\sO(1)^{\oplus p}$ (resp. $\sO(2) \oplus \sO(1)^{\oplus p}$).
            \item[(ii)] Let ${\rm URC}(X)$ be the subset of ${\rm Douady}(X)$ parameterizing all unbendable rational curves on $X$. It is an open subset of ${\rm Douady}(X)$ because  the vector bundle $N_C \cong \sO(1)^{\oplus p} \oplus \sO^{\oplus(n-1-p)}$ on $C \cong \BP^1$ has no nontrivial deformation.  Using  the basic deformation theory of rational curves (e.g. Main Theorem of \cite{Ko}) and $H^1(C, N_C) =0$, we see that ${\rm URC}(X)$ is an open subset in the smooth locus of ${\rm Douady}(X)$. \end{itemize} \end{definition}

                It is convenient to use the following notion.

            \begin{definition}\label{d.jet}
            Let $C$ be a complex manifold and let $\sE$ be a vector bundle on $C$. For a point $x \in C$, denote by $H^0(C, \sE \otimes {\bf m}_x)$ the vector space of sections of $\sE$ vanishing at $x$. The homomorphism $${\rm jet}^{\sE}_x : H^0(C, \sE \otimes {\bf m}_x) \to T^{\vee}_x C \otimes \sE_x$$ is defined by taking the derivative at $x$ of sections vanishing at $x$. \end{definition}

                The following is well known. (See, e.g.,  p.58 of \cite{HM04} and Lemma 3.3 of \cite{HL}).

                \begin{proposition}\label{p.urc}
                In Definition \ref{d.urc}, let $Y$ be a  connected open subset of ${\rm URC}(X).$
                Let $Y \stackrel{\rho}{\leftarrow} Z \stackrel{\mu}{\rightarrow} X$ be the associated universal family morphisms. Define the following distributions  on $Z$
\begin{eqnarray*} \sV &:=& {\rm Ker} ({\rm d} \mu) \\ \sF &:=& {\rm Ker} ({\rm d} \rho) \\
\sT^0 &=& \sV \oplus \sF. \end{eqnarray*}  For  the point $y =[C] \in Y$ corresponding to an unbendable rational curve $C \subset X$ and a point $z \in \rho^{-1}(y) \subset  Z$,  set $\mu(z) = x \in C.$
                   \begin{itemize}
    \item[(i)]  We have the following natural identifications:
          \begin{eqnarray*}T_y Y & = & H^0(C, N_C) \\  T_z Z & = & H^0(C, TX|_C)/H^0(C, TC \otimes {\bf m}_x) \\ \sV_z & = & H^0(C, N_C^+\otimes {\bf m}_x) \\\sF_z & = & H^0(C, TC)/H^0(C, TC \otimes {\bf m}_x) = T_x C, \end{eqnarray*}
           where
           the third (resp. fourth) identification is induced by  the differential ${\rm d}_z \rho: T_z Z \to T_y Y$ (resp. ${\rm d}_z \mu: T_z Z \to T_x X$).
    \item[(ii)] Define $\sT^1:= \p \sT^0.$ It is a vector subbundle of $TZ$ and
the Lie brackets of sections induce a natural isomorphism of vector bundles $\psi: \sF \otimes \sV \to \sT^1/\sT^0$. \item[(iii)] The differential ${\rm d}_z \mu$ sends $\sT_z^1/\sT_z^0$ to $N_{C,x}^+ \subset T_x X/T_xC$ such that in combination with (i) and (ii), we have  the   commutative diagram $$ \begin{array}{ccc}
\sF_z \otimes \sV_z & = &  T_x C \otimes H^0(C, N_C^+ \otimes {\bf m}_x) \\
\psi \downarrow & & \downarrow {\bf j}_x   \\
\sT_z^1/\sT_z^0 & \stackrel{{\rm d}_z \mu}{\longrightarrow} & N_{C,x}^+, \end{array} $$
where ${\bf j}_x$ is the contraction of vectors in $T_x C$ with the image of $$  {\rm jet}^{N_C^+}_x : H^0(C, N_C^+ \otimes {\bf m}_x) \to T^{\vee}_x C \otimes N^+_{C,x}.$$
\item[(iv)] We have isomorphisms $\sF|_{\rho^{-1}(y)} = TC \cong \sO(2) $ and  $$ \sV|_{\rho^{-1}(y)} = \Hom(TC, N_C^+) \cong \sO(-1)^{\oplus m} .$$
\end{itemize}\end{proposition}

            \begin{definition}\label{d.vmrt} In Proposition \ref{p.urc}, the tangent map $\tau: Z \to \BP TX$ sending $z \in Z$ to $[{\rm d} \mu (\sF_z)] \in \BP T_x X$ is an immersion whose image $\sC \subset \BP TX$ is called the {\em variety of minimal rational tangents} (to be abbreviated as VMRT) of the family $Y$.
        The fiber $\sC_x \subset \BP T_x X$ at $x \in \mu(Z)$ is called the {\em VMRT} at $x$. Often, we replace $Y$ by a suitable connected open subset to assume that $\sC \subset \BP TX$ and $\sC_x$ are submanifolds in $\BP TX$.
 \end{definition}

\begin{notation}\label{n.FF}  Let $W$ be a vector space and let $S \subset \BP W$ be a (not necessarily closed) submanifold.  For a point $s\in S$, we denote by $\widehat{s} \subset \widehat{S}$ the corresponding affine cones in $W$. Let  $N_S $ be the normal bundle of $S$ in $\BP W$.
 The second fundamental form
$\II_{S}: \Sym^2 TS \to N_{S}$ is  a homomorphism of vector bundles.  Let ${\rm Dom}(\III_S) \subset S$ be the Zariski-open subset where the image of $\II_S$ is a vector subbundle $N^{(2)}_{S}$ of $N_S$. Then we have the third fundamental form
$$\III_{S,s}: \Sym^3 T_s S \to N_{S,s}/N^{(2)}_{S,s}$$ for each $s \in {\rm Dom}(\III_S).$
Denoting by $\widehat{T}_s S \subset W$ the affine tangent space of $S$ at $s$,
we write $\widehat{N}_{S,s}$ for  $W/\widehat{T}_s S$ and  $\widehat{N}^{(2)}_{S,s}$ for vector subspace of $\widehat{N}_{S,s}$ such that we have natural identifications
$$N_{S,s} = \widehat{s}^{\vee} \otimes \widehat{N}_{S,s} \mbox{ and }
N^{(2)}_{S,s} = \widehat{s}^{\vee} \otimes \widehat{N}^{(2)}_{S,s}.$$
\end{notation}

The following is well-known (see, e.g.,  the proof of Proposition 1.4 in \cite{Hw01},  Corollary 3.14 and Proposition 3.16 in \cite{HL}.)

            \begin{proposition}\label{p.sV}
          In Definitions \ref{d.urc} and \ref{d.vmrt}, replace $Y$ by a neighborhood of $y$ to assume that $\sC$ and $\sC_x$ are submanifolds of $\BP TX$ and the immersion $\tau$ in Definition \ref{d.vmrt} is an embedding.   Let us identify $Z$ with $\sC \subset \BP TX$ by the embedding $\tau$.
\begin{itemize} \item[(i)] The affine tangent space $\widehat{T}_z \sC_x \subset T_x X$
satisfies \begin{eqnarray*} \widehat{T}_z \sC_x & = & {\rm d}_z \mu (\sT^1_z) \\ \widehat{T}_z \sC_x / T_x C & = & N^+_{C,x}.\end{eqnarray*}  \item[(ii)] The natural isomorphism $${\rm d}_z (\tau|_{\mu^{-1}(x)}): \sV_z = H^0(C, N^+_C \otimes {\bf m}_x)  \ \to \ T_z \sC_x = T^{\vee}_x C \otimes N^+_{C,x}$$ coincides with ${\rm jet}^{N_C^+}_x$.
\item[(iii)] For  a local section $\vec{f}$ of $\sF$ and local sections $\vec{v}_1, \vec{v}_2$ of $\sV$ near $z \in \sC$ with  values $v_1, v_2 \in \sV_z$ at $z$, $$ ([\vec{v}_2, [\vec{v}_1, \vec{f}]]_z  \mod \sT^1_z)  \ = \   \vec{f}_z \otimes \II_{\sC_x, z}(v_1, v_2) \  $$
where the left hand side is regarded as an element of $\widehat{N}_{\sC_x, z} = N_{C,x}/N_{C,x}^+$ by (i) and  Proposition \ref{p.urc} (iii), and $\vec{f}_z \otimes$ stands for the contraction of $\vec{f}_z \in \sF_z = T_x C$ with $$\II_{\sC_x, z}(v_1, v_2) \in T^{\vee}_x C \otimes \widehat{N}^{(2)}_{\sC_x, z} \subset T^{\vee}_x C \otimes \widehat{N}_{\sC_x,z}.$$
\item[(iv)] Assume that $z \in \sC_x$ is in ${\rm Dom}(\III_{\sC_x})$. Then $\sT^2 = \p \sT^1 = \p^{(2)}\sT^0$ is a vector subbundle of $TZ$ in a neighborhood of $\rho^{-1}(y)$ in $Z$ and
$$\widehat{N}^{(2)}_{\sC_x, z}  = ( {\rm d}_z \mu (\sT^2_z) \mod \widehat{T}_z \sC_x)$$
as subspaces of $\widehat{N}_{\sC_x, z} = N_{C,x}/N^+_{C,x}.$
    \item[(v)] In (iv),  for  a local section $\vec{f}$ of $\sF$ and local sections $\vec{v}_1, \vec{v}_2, \vec{v}_3$ of $\sV$ near $z \in \sC$ with  values $v_1, v_2, v_3 \in \sV_z$ at $z$,  $$ ([\vec{v}_3, [\vec{v}_2, [\vec{v}_1, \vec{f}]]]_z \mod \sT^2_z) =   \vec{f}_z \otimes \III_{\sC_x, z}(v_1, v_2, v_3) , $$
where the left hand side is regarded as an element of $\widehat{N}_{\sC_x,z}/\widehat{N}^{(2)}_{\sC_x,z}$ by (iv), and $\vec{f}_z \otimes$ stands for the contraction of $\vec{f}_z \in \sF_z = T_x C$ with $$ \III_{\sC_x, z}(v_1, v_2, v_3) \in  T^{\vee}_x C \otimes \widehat{N}_{\sC_x,z}/\widehat{N}^{(2)}_{\sC_x, z}.$$ \end{itemize}\end{proposition}

\begin{definition}\label{d.tensor}
In Proposition \ref{p.urc},
define $\sD_y := H^0(C, N_C^+) \subset T_y Y$. This determines a distribution $\sD \subset TY$ of rank $2m.$
For each $x \in C,$ define $$U_x := H^0(C, \sO(x)),$$ the 2-dimensional vector space of rational functions on $C$ with at most one pole at $x$ and define $$V_x := H^0(C, N_C^+\otimes {\bf m}_x) \subset \sD_y,$$ such that we have a canonical tensor decomposition $\sD_y = U_x \otimes V_x$.
\end{definition}

\begin{proposition}\label{p.T1D}
In Definition \ref{d.tensor}, let $\mathbf{1}_x \in U_x$ be the constant function on $C$ with value 1. The differential ${\rm d}_z \rho: T_z Z \to T_y Y$  induces identification $\sT^1_z/\sF_z  = \sD_y$  such that
 $$\sV_z = \mathbf{1}_x \otimes  V_x \ \subset \ U_x \otimes V_x = \sD_y.$$
 Moreover, when $ z\in {\rm Dom}(\III_{\sC_x}),$ it induces an identification $\sT^2_z/\sF_z = (\p \sD)_y.$
\end{proposition}

\begin{proof}
The identification $\sT^1_z/\sF_z  = \sD_y$  follows from Proposition 1 and Proposition 8 of \cite{HM04} (see Proposition 3.7 of \cite{HL}). Then $\sV_z = \mathbf{1}_x \otimes  V_x$ follows from Proposition \ref{p.urc} (i) and the identification $\sT^2_z/\sF_z = (\p \sD)_y$ comes from $\sT^2 = \p \sT^1$ in Proposition \ref{p.sV} (iv).  \end{proof}

We skip the proof of the following two elementary lemmata.

\begin{lemma}\label{l.res}
In Proposition \ref{p.T1D}, let $f \in U_x$ be a nonconstant rational function on $C$. \begin{itemize} \item[(i)] The  homomorphism $${\bf res}_x: \wedge^2 H^0(C, \sO(x)) \to T_x C$$ that sends ${\bf 1}_x \wedge f$ to $${\rm Res}_x ( f \ {\rm d}t ) \frac{\p}{\p t},$$
where $t$ is a local holomorphic coordinate on $C$ centered at $x$ and ${\rm Res}_x( f \ {\rm d} t)$ is
the residue of the logarithmic form $f \ {\rm d}t$ at $x$, is independent of the choice of the coordinate $t$ and gives a canonical identification of $\wedge^2 U_x$ with $T_x C$.
\item[(ii)] Write  $\vec{f}_x := {\bf res}_x ({\bf 1}_x \wedge f) \in T_x C$. For any vector bundle $\sE$ on $C$ and a section $v \in H^0(C, \sE \otimes {\bf m}_x)$,  let $ \vec{f}_x (v) \in \sE_x$ be the contraction of  $\vec{f}_x$ with $ {\rm jet}^{\sE}_x(v) \in T_x^{\vee} C \otimes \sE_x.$ Then $$ \vec{f}_x (v) = (f \otimes v) (x) \ \in \sE_x $$ where $f \otimes v $ stands for the section of $H^0(C, \sE)$  under the natural homomorphism $$U_x \otimes H^0(C, \sE \otimes {\bf m}_x) \to H^0(C, \sE).$$
\end{itemize} \end{lemma}

           \begin{lemma}\label{l.tensor}
In Lemma \ref{l.res}, for any point $x' \in C$ different from $x$, define $$A_x^{x'} := \{ f \in U_x \mid f(x') = 0\} \subset U_x.$$ Then we have canonical isomorphisms $$A_{x}^{x'} \otimes U_{x'} = U_x \mbox{ and } A_x^{x'} \otimes V_x = V_{x'}$$ coming from multiplication by rational functions.  Moreover,  for any $f \in A_x^{x'}$ and $v \in V_x$, we have $f \otimes v = {\bf 1}_{x'} \otimes f v$ under the two tensor decompositions $\sD_y = U_x \otimes V_x = U_{x'} \otimes V_{x'}.$  \end{lemma}

\begin{proposition}\label{p.f}
In Definition \ref{d.tensor}, for a local section $\vec{f}$ of $\sF$ and $\vec{v}$ of $\sV$ near $z \in \sC$, denote by  $f \in U_x$ a rational function satisfying ${\rm d}_z \mu (\vec{f}_z) = {\bf res}_x ({\bf 1}_x \wedge f)$ in Lemma \ref{l.res} and     by $v \in \sV_z = V_x$ the value of $\vec{v}$ at $z$. Then \begin{equation}\label{e.f} {\rm d}_z \rho ([\vec{v}, \vec{f}]_z) \ \equiv \  f \otimes v   \mod {\mathbf 1}_x \otimes V_x ( \subset U_x \otimes V_x = \sD_y). \end{equation} \end{proposition}

\begin{proof}
The vector $[\vec{v}, \vec{f}]_z \in T_z Z$ is represented by a section $w \in H^0(C, TX)$ by Proposition \ref{p.urc} (i). Regarding $v$ as an element of $H^0(C, N^+_C \otimes \mathbf{m}_x)$, Proposition \ref{p.urc} (iii) says that $$({\rm d}_z \mu ([\vec{v}, \vec{f}]_z) \mod T_xC) = \vec{f}_z \otimes {\rm jet}^{N_C^+}_x(v) = \vec{f}_z (v) \in N^+_{C,x}.$$ This implies that the value $w_x$ at $x$ of the section $w$ of $TX|_C$ satisfies
 $$(w_x  \mod T_xC) = \vec{f}_z (v).$$ Since ${\rm d}_z \rho([\vec{v}, \vec{f}]_z)$ in $T_y Y = H^0(C, N_C)$ is represented by $w$ modulo $H^0(C, TC)$, its value modulo ${\bf 1}_x \otimes V_x = H^0(C, N_C \otimes {\bf m}_x)$ is just $$ (w_x \mod T_x C) \in H^0(C, N_C)/H^0(C, N_C \otimes {\bf m}_x) = N_{C,x}.$$ So the left hand side of (\ref{e.f}) is just $\vec{f}_z(v)$ in Lemma \ref{l.res} (ii). Then Lemma \ref{l.res} (ii) says that this is equal to $(f\otimes v) (x)$,  the right hand side of (\ref{e.f}).
\end{proof}
%

\begin{proposition}\label{p.II}
In Proposition \ref{p.f}, for $u_1, u_2 \in U_x$ and $v_1, v_2 \in V_x$, let $u_1 \otimes v_1, u_2 \otimes v_2$ be elements of $\sD_y = U_x \otimes V_x$. Then we have \begin{equation}\label{e.Levitensor} {\rm Levi}^{\sD}_y(u_1 \otimes v_1, u_2 \otimes v_2) = (u_1 \wedge u_2) \otimes \II_{\sC_x, z} (v_1, v_2), \end{equation} where the second fundamental form $\II_{\sC_x,z}: \Sym^2 V_x \to T^{\vee}_x C \otimes (N_{C,x}/N^+_{C,x})$ is interpreted as a homomorphism $$\Sym^2 V_x \to \wedge^2 U_x^{\vee} \otimes (H^0(C, N_C)/H^0(C, N_C^+)) = \wedge^2 U_x^{\vee} \otimes (T_y Y/\sD_y)$$ via the isomorphism ${\bf res}_x: T_x C \cong \wedge^2 U_x$ in Lemma \ref{l.res} (i) and the natural isomorphism $
H^0(C, N_C/N_C^+) = N_{C,x}/N_{C,x}^+$ coming from $N_C/N_C^+ \cong \sO^{\oplus (n-m-1)}$. \end{proposition}

\begin{proof}
We may check  (\ref{e.Levitensor}) assuming $u_1 = f$ and $u_2= {\bf 1}_x$. It is a direct consequence of
\begin{equation}\label{e.2.3} ([\vec{v}_2, [\vec{v}_1, \vec{f}]]_z  \mod \sV_z)  \ = \   \vec{f}_z \otimes \II_{\sC_x, z}(v_1, v_2)   \end{equation}
from Proposition \ref{p.sV} (iii). In fact,        the vector $v_2 \in \sV_z$ is sent by ${\rm d}_z \rho$ to the corresponding element $v_2 \in V_x = H^0(C, N_C^+ \otimes {\bf m}_x)$, which is just ${\bf 1}_x \otimes v_2$ in $U_x \otimes V_x = \sD_y$.   The vector $ [\vec{v}_1, \vec{f}]_z $ is sent to $f \otimes v_1$ modulo $\mathbf{1}_x \otimes V_x$ by Proposition \ref{p.f}. Then the left hand side of (\ref{e.2.3}) is equal to ${\rm Levi}^{\sD}(u_2 \otimes v_2, u_1 \otimes v_1)$.  By Lemma \ref{l.res}, the isomorphism $\mathbf{res}_x$ identifies    $\vec{f}_z$ with  $-\mathbf{1}_x \wedge f = - u_1 \wedge u_2$. Thus the right hand side of (\ref{e.2.3}) is  $$- (u_1 \wedge u_2) \otimes \II_{\sC_x, z}(v_1, v_2).$$ This proves (\ref{e.Levitensor}).  \end{proof}

\begin{proposition}\label{p.III}
In Proposition \ref{p.f}, assume that $\sD$ is regular at $y \in Y$. Then for any $f \in U_x, v_1, v_2, v_3 \in V_x$,  the Lie bracket in ${\rm symb}_y  (\sD)$ satisfies $$[ \mathbf{1}_x \otimes v_3, [ \mathbf{1}_x \otimes v_2, f \otimes v_1] ]= (\mathbf{1}_x \wedge f) \otimes \III_{\sC_x, z}(v_1, v_2, v_3),$$
where  $\III_{\sC_x, z}: \Sym^3 V_x \to T^{\vee}_x C \otimes (\widehat{N}_{\sC_x,z}/\widehat{N}^{(2)}_{\sC_x,z})$ is interpreted as a homomorphism $$\Sym^3 V_x \to \wedge^2 U^{\vee}_x \otimes (T_y Y/(\p \sD)_y),$$ via   the isomorphism ${\bf res}_x: T_x C \cong \wedge^2 U_x$ and the isomorphisms $$\widehat{N}_{\sC_x,z}/\widehat{N}^{(2)}_{\sC_x,z} \stackrel{{\rm d}_z \mu}{= }   T_z z/\sT^2_z \stackrel{{\rm d}_z \rho}{=} T_y Y/(\p \sD)_y.$$ \end{proposition}

\begin{proof}
The proof follows  the same argument as  the proof of Proposition \ref{p.II}, using Proposition \ref{p.sV} (v) in place of Proposition \ref{p.sV} (iii). \end{proof}
%

\section{From nondegenerate lines on contact manifolds to distributions with symbols $\fg_+(F)$}\label{s.contact}

\begin{definition}\label{d.line}
Let $X$ be a complex manifold of dimension $2m +3, m \geq 1,$ and let $\sH \subset TX$ be a contact distribution, namely, a distribution of rank $2m+2$ such that the Levi tensor ${\rm Levi}^{\sH}: \wedge^2 \sH \to TX/\sH$ gives a symplectic form, i.e., a nondegenerate anti-symmetric form, on each fiber of $\sH_x$. The pair $(X, \sH)$ is called a {\em contact manifold}. The quotient line bundle $\sL:= TX/\sH$ is called the {\em contact line bundle} on $X$.   An unbendable smooth rational curve $C \subset X$ satisfying $\sL|_C \cong \sO(1)$ is called a {\em  line}.   Denote by ${\rm Lines}(X, \sH)$  the open subset of ${\rm URC}(X) \subset {\rm Douady}(X)$ parameterizing lines. \end{definition}

\begin{definition}\label{d.Legendrian}
Fix a 1-dimensional vector space $\mathbf{L}$ and let $\omega: \wedge^2 W \to \mathbf{L}$ be a symplectic form on a vector space $W$ of dimension $2m+2$.
\begin{itemize}
\item[(i)] For a subspace $B \subset W$, define $B^{\perp} := \{ w \in W \mid \omega(w, B) =0\}.$ \item[(ii)] A subspace $B \subset W$ is {\em Lagrangian} if $B^{\perp} = B$.
    In this case, the dimension of $B$ is $m+1$.
\item[(iii)]
A submanifold $S \subset \BP W$ is {\em Legendrian} if its affine tangent space $\widehat{T}_s S \subset W$ at every point $s \in S$ is a Lagrangian subspace of $W$.
\end{itemize}
\end{definition}

\begin{lemma}\label{l.Legendrianvmrt}
In Definition \ref{d.line}, choose a connected open subset $Y \subset {\rm Lines}(X, \sH)$ such that $\sC \subset \BP TX$ in Definition \ref{d.vmrt} is a submanifold.
\begin{itemize} \item[(i)] For each line $C \subset X$, the normal bundle $N_C$ is isomorphic to $\sO(1)^{\oplus m} \oplus \sO^{\oplus (m+2)}$ and $\sH|_C$ is isomorphic to
$\sO(2) \oplus \sO(1)^{\oplus m} \oplus \sO^{\oplus m} \oplus \sO(-1)$. In other words, the integer $p$ in Definition \ref{d.urc} is exactly $m = \frac{1}{2}(\dim X-3)$.
\item[(ii)] For each $x \in C$, the VMRT $\sC_x \subset \BP T_x X$ is contained in  $\BP \sH_x$ and is Legendrian with respect to the symplectic form  ${\rm Levi}^{\sH}_x$ on $\sH_x$.
    \end{itemize} \end{lemma}

\begin{proof}
From the well-known relation $\det TX = \sL^{\otimes (m+2)}$ (for example, from (2.2) of \cite{LB}), we have the isomorphism $N_C \cong \sO(1)^{\oplus m} \oplus \sO^{\oplus (m+2)}.$ Consequently, $$TX|_C \cong TC \oplus N_C \cong  \sO(2) \oplus \sO(1)^{\oplus m} \oplus \sO^{\oplus (m+2)}.$$  Since $TC \cong \sO(2)$, it must belong to the kernel of the projection $TX|_C \to (\sL = TX/\sH)|_C \cong \sO(1).$ Since this holds for all lines represented by elements in the parameter space $Y$, we have $\sC\subset\BP\sH$. By Proposition \ref{p.sV}(i),  $\widehat{T}_z\mathcal{C}_x/T_xC=N^+_{C, x}$, where $x\in C$ and $z=\mathbb{P}T_xC\in\mathcal{C}_x$. It follows that $\sO(2) + \sO(1)^{\oplus m}$ is contained in $\sH|_C$. This gives $\sH|_C \cong \sO(2) \oplus \sO(1)^{\oplus m} \oplus \sO^{\oplus m} \oplus \sO(-1)$, completing the proof of (i).
By (i),  any line is tangent to $\sH$, and $\sC_x \subset \BP \sH_x$. Moreover, the Levi tensor gives a family of nondegenerate antisymmetric forms on the bundle $\sH$, i.e., a homomorphism  $$\wedge^2 \sH_C \cong \wedge^2(\sO(2) \oplus \sO(1)^{\oplus m} \oplus \sO^{\oplus m} \oplus \sO(-1)) \to \sL|_C \cong \sO(1).$$
The subspace $(\sO(2) \oplus \sO(1)^{\oplus m})_x$ is a Lagrangian subspace of $\sH_x$ with respect to this antisymmetric form. As this subspace corresponds to the affine tangent space of $\sC_x$ at $\BP T_x C$ by Proposition \ref{p.sV} (i), the VMRT must be Legendrian, proving  (ii). \end{proof}

\begin{definition}\label{d.sR}
In Definition \ref{d.line}, let $Y \subset {\rm Lines}(X, \sH)$ be a connected open subset
and $Y \stackrel{\rho}{\leftarrow} Z \stackrel{\mu}{\rightarrow} X$ be the universal family.
Given $z \in Z$, write $x= \mu(z), y = \rho(z)$ and  $ C = \mu(\rho^{-1}(y)) \subset X$.   Using the symplectic form ${\rm Levi}_x^{\sH}$ on $\sH_x$, define $$\sR_z := (T_xC)^{\perp} \ \subset \sH_x \mbox{ and } \sR^+_z : = \widehat{T}_z \sC_x  \ \subset \sR_z \subset \sH_x.$$
Then $\sR^+ \subset \sR $ are vector subbundles of the vector bundle  $\sH$ on $Z$,  of rank $m+1$ and $2m+1$, respectively.
Denote by $\sR^+_C \subset \sR_C \subset \sH|_C$ the corresponding vector subbundles on $C$.  \end{definition}

\begin{lemma}\label{l.secH}
In Definition \ref{d.sR}, \begin{itemize} \item[(i)] the vector bundle $\sR_C$ (resp. $\sR^+_C$)  is isomorphic to the subbundle
$$\sO(2) \oplus \sO(1)^{\oplus m} \oplus \sO^{\oplus m} \ (\mbox{ resp. } \sO(2) \oplus \sO(1)^{\oplus m})$$ of $\sH|_C$ from Lemma \ref{l.Legendrianvmrt} (i); \item[(ii)] if $s \in H^0(C, TX|_C)$ satisfies $s_x \in \sR_{C,x}$ for a point $x \in C$, then $s \in H^0(C, \sH|_C);$ and
\item[(iii)]  ${\rm Levi}^{\sH}$ induces a perfect pairing of vector bundles on $C$
    $$ N_C^+ \otimes \sR_C/(TX|_C)^+  \to \sL|_C.$$  \end{itemize} \end{lemma}

\begin{proof}
Note that ${\rm Levi}^{\sH}$ gives a perfect paring $$ TC \otimes \sH|_C/\sR_C  \cong \sL|_C \cong \sO(1).$$ This implies (i).
From (i), the quotient bundle $TX|_C/\sR_C$ is isomorphic to $\sO^{\oplus 2}$. The assumption $s_x \in \sR_{C,x}$ implies that $s$ modulo $\sR_C$ defines a section of $\sO^{\oplus 2}$ vanishing at $x$. Hence it vanishes identically. This shows that $s \in H^0(C, \sH) = H^0(C, \sR)$, proving (ii). (iii) is immediate from (i).
\end{proof}

\begin{definition}\label{d.nondegenerate}
A cubic form $F \in \Sym^3 V^{\vee}$ on a vector space $V$ is {\em nondegenerate}
if $$\{ v \in V \mid F(v, v_1, v_2) = 0 \mbox{ for all } v_1, v_2 \in V\} = 0.$$
\end{definition}

\begin{definition}\label{d.LM}
In Definition \ref{d.Legendrian}, let $S \subset \BP W$ be a Legendrian submanifold.
We say that $S$ has {\em nondegenerate fundamental forms} at $s \in S$ if
\begin{itemize} \item[(i)]  the point $s$ is contained in $ {\rm Dom}(\III_S)$;\item[(ii)] the image of the second fundamental form $\widehat{N}^{(2)}_{S,s}$ is equal to  $\widehat{s}^{\perp}/\widehat{T}_s S$;  \item[(iii)] the third fundamental form $\III_{S,s}: \Sym^3 T_s S \to \widehat{s}^{\vee} \otimes W/\widehat{s}^{\perp}$ is a nondegenerate cubic form on $T_s S$ in the sense of Definition \ref{d.nondegenerate}; and \item[(iv)]
for any $v_1, v_2 \in T_s S$, the element $$\III_{S,s}(v_1, v_2, \cdot) \in T^{\vee}_sS \otimes \widehat{s}^{\vee} \otimes W/\widehat{s}^{\perp}$$ coincides with $$\II_{S,s}(v_1, v_2) \in \widehat{s}^{\vee} \otimes \widehat{N}^{(2)}_{S,s} = \widehat{s}^{\vee}\otimes \widehat{s}^{\perp}/\widehat{T}_s S$$ via the natural isomorphisms $$\widehat{s}^{\perp}/\widehat{T}_s S = (\widehat{T}_s S/\widehat{s})^{\vee} \otimes \mathbf{L} \mbox{ and } W/\widehat{s}^{\perp} = \widehat{s}^{\vee} \otimes \mathbf{L}$$ induced by $\omega$.    \end{itemize} \end{definition}

\begin{definition}\label{d.Gauss}
A line $C$ on a contact manifold $(X, \sH)$ is a {\em nondegenerate line},  if for some point $x\in C$, the VMRT  $\sC_x$ at $x$ has nondegenerate fundamental forms at the point $\BP T_x C \in \sC_x$. Let ${\rm NDL}(X, \sH)$ be the subset of ${\rm Lines}(X, \sH)$ parameterizing nondegenerate lines. It is a Zariski-open subset in ${\rm Lines}(X, \sH)$.  \end{definition}

\begin{remark}\label{r.LM}
The conditions in Definition \ref{d.LM} may look technical, but actually they hold at general points of a general Legendrian submanifold. For example, Section 3 of \cite{LM} shows that
if the Legendrian submanifold $S \subset \BP W$ is a smooth projective variety different from a linear subspace, then $S$ has nondegenerate fundamental forms at a general point $s \in S$.
When $(X, \sH)$ is a Fano contact manifold of Picard number 1, Theorem 1.1 of \cite{Ke2} implies that the VMRT at a general point is a smooth projective variety.  Thus  a general line is nondegenerate in the sense of Definition \ref{d.Gauss}, unless the VMRT's are linear. It is expected that the latter situation does not occur and \cite{BKK} proposes an approach to exclude this possibility. \end{remark}

\begin{proposition}\label{p.ndII}
In Definitions  \ref{d.Gauss}, assume that $Y \subset {\rm NDL}(X, \sH).$ Then
the second fundamental forms of the VMRT determine a surjective homomorphism $\II: \Sym^2 \sV \to \sF^{\vee} \otimes \sR/\sR^+$  and the third fundamental forms of the VMRT determine a surjective homomorphism $\III: \Sym^3 \sV \to \sF^{\vee} \otimes (\sH/\sR).$
In particular, we have $$\sT^1_z = ({\rm d}_z \mu)^{-1}(\sR^+_z) \subset T_z Z \mbox{ and } \sT^2_z = ({\rm d}_z \mu)^{-1} (\sR_z) \subset T_z Z$$ for any $z \in Z$. \end{proposition}

\begin{proof}
It is immediate from the definition of a nondegenerate line that the second and the third fundamental forms of VMRT determine homomorphisms $\II$ and $\III,$  which are surjective at general points of $Z$.  By Proposition \ref{p.urc} (iv) and
$ (\sH|_C/\sR_C) \cong \sO(-1)$ from Lemma \ref{l.secH} (i), the two homomorphisms are surjective at every point of $Z$. Then $\sT^2_z = ({\rm d}_z \mu)^{-1} (\sR_z)$ follows from Proposition \ref{p.sV} (iv). \end{proof}

\begin{proposition}\label{p.psD}
Under the hypotheses of Proposition \ref{p.ndII},  the two subspaces $(\p \sD)_y$ and $H^0(C, \sH|_C/TC)$ of $T_y Y = H^0(C, N_C)$ coincide.
Consequently, the two quotient spaces
$T_y Y /(\p \sD)_y$ and  $$H^0(C, \sL) = H^0(C, N_C)/H^0(C, \sH|_C/TC)$$ can be identified in a natural way.
\end{proposition}

\begin{proof}
By Lemma \ref{l.secH} (ii), the subspace $\sT_z^2 = ({\rm d}_z \mu)^{-1} \sR_z \subset T_z Z$ in Proposition \ref{p.ndII} can be identified with $$ H^0(C, \sH)/H^0(C, TC \otimes {\bf m}_x) \ \subset \ H^0(C, TX)/H^0(C, TC \otimes {\bf m}_x)$$ under the identification in Proposition \ref{p.urc} (i).   Since ${\rm d} \rho: T_z Z \to T_y Y$ is the quotient map $$H^0(C, TX)/H^0(C, TC \otimes {\bf m}_x) \to H^0(C, N_C),$$ the subspace $({\rm d}_z \mu)^{-1} \sR$ is sent to
$H^0(C, \sH|_C/TC) \subset T_y Y$. This agrees with $(\p \sD)_y$  because ${\rm d}_z \rho (\sT^2_z) = \p\sD$ by Proposition \ref{p.T1D}.
\end{proof}


\begin{lemma}\label{l.parallelIII}
In Proposition \ref{p.psD},
 define $$I_x := T_x^{\vee}C \otimes H^0(C, \sL \otimes {\bf m}_x) = (T_x^{\vee}C)^{\otimes 2} \otimes \sL_x$$ and  write $\mathbf{F}_x: \Sym^3 V_x \to I_x$
for $\III_{\sC_x, z}$.
  For each pair $x \neq x' \in C$, recall from Lemma \ref{l.tensor} the 1-dimensional vector space $A_x^{x'}$ of rational functions on $C$ with pole at $x$ and zero at $x'$.
 \begin{itemize} \item[(i)]
There is a canonical isomorphism between $(A^{x'}_x)^{\otimes 3} \otimes I_x$ and $I_{x'}$. For $f \in A^{x'}_x$ and $j \in I_x$, denote by $f^3 \cdot j$ the  element of $I_{x'}$ corresponding to $ f^3 \otimes j  \in (A^{x'}_x)^{\otimes 3} \otimes I_x$. \item[(ii)]
For any $v_1, v_2, v_3 \in V_{x}$ and $f \in A^{x'}_x$,
$$f^3 \cdot \mathbf{F}_x(v_1, v_2, v_3)  =  \mathbf{F}_{x'}(f \cdot v_1, f \cdot v_2, f \cdot  v_3).$$
\end{itemize}
\end{lemma}

\begin{proof}
Fix a base point $x_0 \in C$ and define a line bundle $\sA$ on $C$ whose fiber at $x$ is $A_x^{x_0}$. Then $\sA$ is isomorphic to $\sO(1)$.
 Let  $\sI$ be the line bundle on $C$ whose fiber at $x \in C$ is  $I_x$. Then $\sI$ is isomorphic to $\sO(-3)$.
  Thus  we have canonical trivializations of the line bundle $\sA^{\otimes 3} \otimes \sI = I_{x_0} \times C,$ which  shows (i).

Note that $\sV|_{\rho^{-1}(y)} \cong \sO(-1)^{\oplus m}$ from Proposition \ref{p.urc} (iv).
Thus the collection of the homomorphisms $\mathbf{F}_x$ for all $x\in C$  gives rise to a surjective homomorphism between trivial vector bundles $ \Sym^3 (\sV \otimes \sA) \to \sA^{\otimes 3} \otimes \sI,$ which corresponds to $\mathbf{F}_{x_0}: \Sym^3 V_{x_0} \to I_{x_0}$ under the canonical trivialization $\sA^{\otimes 3} \otimes \sI$. This shows (ii). \end{proof}

\begin{proposition}\label{p.sD}
In Proposition \ref{p.psD},    there are natural identifications
\begin{itemize} \item[(i)] $\sD_y  =  U_x \otimes V_x; $
\item[(ii)] $(\p \sD)_y/\sD_y  =  \wedge^2 U_x \otimes I_x \otimes V^{\vee}_x,$ and
\item[(iii)] $T_y Y/ (\p \sD)_y  = \wedge^2 U_x \otimes I_x \otimes U_x.$ \end{itemize}
In particular, for any pair $x \neq x' \in C$, there is a natural identification of $ \wedge^2 U_x \otimes I_x \otimes U_x$ and $\wedge^2 U_{x'} \otimes I_{x'} \otimes U_{x'}$
compatible with the identification $(\sA_{x}^{x'})^{\otimes 3} \otimes I_x = I_{x'}$ in Lemma \ref{l.parallelIII} and the identification $\sA_{x'}^x \otimes U_{x'} = U_x$   in Lemma \ref{l.tensor}.   \end{proposition}

\begin{proof}
We have already seen (i) in Definition \ref{d.tensor}. (ii) is from $(\p \sD)_y/\sD_y = \sT_z^2/\sT_z^1$ in Proposition \ref{p.T1D},
 the identification $\wedge^2 U_x = T_x C$ from Lemma \ref{l.res}, and
$$\sT_z^2/\sT_z^1 = (N^+_{C,x})^{\vee} \otimes \sL_x = V^{\vee}_x \otimes T^{\vee}_x C \otimes \sL_x $$ from Lemma \ref{l.secH} (iii).  (iii) follows from $$T_yY/(\p \sD)_y = H^0(C, \sL) = H^0(C, \sL \otimes {\bf m}_x) \otimes U_x,$$
where the first equality is by Proposition \ref{p.psD}. \end{proof}

\begin{definition}\label{d.fg+}
Let $\U$ be a vector space of dimension 2 and let $\V$ be a vector space of dimension $m \geq 1.$
Fix a 1-dimensional vector space $\I$.
Let $\fg_+ = \fg_1 \oplus \fg_2 \oplus \fg_3$ be the vector space defined by  $$\fg_1 := \U \otimes \V, \ \fg_2 := (\wedge^2 \U) \otimes \I \otimes  \V^{\vee}, \ \fg_3 := (\wedge^2 \U) \otimes \I  \otimes \U.$$
Fix a  cubic form $F: \Sym^3 \V \to \I$ on $\V$ which is nondegenerate in the sense of Definition \ref{d.nondegenerate}.
For $v_1, v_2 \in \V$, define $F_{v_1 v_2} \in \I \otimes \V^{\vee}$ by $$ F_{v_1 v_2} (v) : = F(v_1, v_2, v) \in \I.$$
We define a graded Lie algebra structure on
 $\fg_+ := \fg_1 \oplus \fg_2 \oplus \fg_3$ by
    \begin{equation}\label{e.fg1} [u_1 \otimes v_1, u_2 \otimes v_2] = (u_1 \wedge u_2) \otimes F_{v_1, v_2} \mbox{ and } \end{equation}
\begin{equation}\label{e.fg2} [[u_1 \otimes v_1, u_2 \otimes v_2], u_3 \otimes v_3] = (u_1 \wedge u_2) \otimes  F(v_1, v_2, v_3)  \otimes u_3.\end{equation}
By the nondegeneracy of $F$, (\ref{e.fg1}) implies $[\fg_1, \fg_1] = \fg_2$. Thus (\ref{e.fg2}) is sufficient to determine the Lie bracket $[\fg_1, \fg_2]$. It is easy to check that this gives a graded Lie algebra structure on $\fg_+$. Note that the Lie bracket $[\fg_1, \fg_2] \to \fg_3$ is independent of $F$ and satisfies
      \begin{equation}\label{e.fg3}  [(u_1 \wedge u_2) \otimes v^*, u_3 \otimes v_3]  =  (u_1\wedge u_2) \otimes u_3 \otimes v^*(v_3) \end{equation}
     for any $u_1, u_2, u_3 \in \U,  v_3 \in \V$ and $v^* \in \I \otimes V^{\vee}$.
    Sometimes, we write $\fg_+ = \fg_+(F)$ to indicate that the Lie algebra structure depends on the cubic form $F$.
    \end{definition}

The following is the precise formulation of going from the right hand side to the left hand side of (\textsf{Corr.$m$}) in Section \ref{s.1}.

\begin{theorem}\label{t.symbol}
Let $(X, \sH)$ be a contact manifold of dimension $2m+3$.
Let $C \subset (X, \sH)$ be a nondegenerate  line  and let $y \in {\rm NDL}(X, \sH)$ be the corresponding point. Then
the symbol algebra ${\rm symb}_y (\sD)$ of the natural distribution $\sD$ on ${\rm NDL}(X, \sH)$ is isomorphic to $\fg_+(F_y)$ in Definition \ref{d.fg+} for some nondegenerate cubic form $F_y$ on an $m$-dimensional vector space $\V$.
\end{theorem}

\begin{proof}
Let $\mathbf{F}_x: \Sym^3 V_x \to I_x$ be the cubic form given by $$\III_{\sC_x, z}: \Sym^3 \sV_x \to (\sF_z^{\vee})^{\otimes 2} \otimes \sL_x$$ via the natural isomorphism $V_x \cong \sV_x$ for a point  $x \in C$.  Using
Proposition \ref{p.sD}, let us identify ${\rm symb}_y(\sD)$ with $\fg_+$ as graded vector spaces by setting $\U = U_{x}, \V = V_{x}$ and $\I = I_{x}$.
It suffices to show that the Lie algebra structure of ${\rm symb}_y(\sD)$ agrees with  that of $\fg_+(\mathbf{F}_{x}),$ namely, it satisfies (\ref{e.fg1}) and (\ref{e.fg2}).

For any $x \in C$ and   $v_1, v_2 \in V_x$,
let $F_{v_1 v_2} \in I_x \otimes V_x^{\vee}$ be  the contraction of $F= \mathbf{F}_{x}$ with $v_1, v_2$. Fix a nonconstant rational function $f \in  U_x$.    For $v_1,  v_2 \in V_{x},$ Proposition \ref{p.II}, combined with Definition \ref{d.LM} (iv),  says that in ${\rm symb}_y(\sD)$,
$$ [{\bf 1}_x \otimes v_2, f \otimes v_1] = ({\bf 1}_x \wedge f) \otimes F_{v_1 v_2}  \ \in (\wedge^2 U_x) \otimes (I_x \otimes V^{\vee}_x). $$ Thus the Lie algebra ${\rm symb}_y(\sD)$ satisfies (\ref{e.fg1})  with $F= \mathbf{F}_{x}$.

For $v_1, v_2, v_3 \in V_x,$ Proposition \ref{p.III} says
\begin{equation}\label{e.p2}  [{\bf 1}_x \otimes v_3, [ {\bf 1}_x \otimes v_2, f \otimes v_1]]  =  ({\bf 1}_x \wedge f) \otimes \mathbf{F}_x(v_1, v_2, v_3) \otimes {\bf 1}_x \end{equation} as elements in $ (\wedge^2 U_x) \otimes I_x \otimes U_x. $ Thus the Lie algebra ${\rm symb}_y(\sD)$ satisfies  (\ref{e.fg2}) with $F= \mathbf{F}_{x}$ when $  u_3 = {\bf 1}_{x}$. It remains to check (\ref{e.fg2}) when $u_3 $ is a nonconstant function in $ U_{x}.$

For a given nonconstant function $f \in U_{x}$, let $x' \in C$ be the zero of $f$ such that $f \in \sA^{x'}_x$ in the notation of Lemma \ref{l.parallelIII}.   For $v_1, v_2, v_3 \in V_{x}$, we have  $ w_1, w_2, w_3  \in  V_{x'}$  such that $w_i = f \cdot v_i$ for $i=1,2,3$. Let $h \in \sA^x_{x'}$ be the rational function $\frac{1}{f}$.   Then \begin{equation}\label{e.1tensor} {\bf 1}_x \otimes v_i  = h \otimes w_i \mbox{ and }  {\bf 1}_{x'} \otimes w_i  = f \otimes v_i     \end{equation} for $i=1,2,3,$ under the tensor decomposition $H^0(C, N_C^+) = U_x \otimes V_x = U_{x'} \otimes V_{x'}$.
From (\ref{e.p2}) applied to the point $x' \in C$,
\begin{equation}\label{e.2tensor}  [{\bf 1}_{x'} \otimes w_3, [ {\bf 1}_{x'} \otimes w_2, h \otimes w_1]]  =  ({\bf 1}_{x'} \wedge h) \otimes \mathbf{F}_{x'}(w_1, w_2, w_3) \otimes {\bf 1}_{x'} \end{equation} as elements of $T_y Y/(\p \sD)_y =  (\wedge^2 U_{x'}) \otimes I_{x'} \otimes U_{x'}.$
By (\ref{e.1tensor}), the left hand side of (\ref{e.2tensor}) is equal to $  [f \otimes v_3, [f \otimes v_2, {\bf 1}_x \otimes v_1] ].$ On the other hand, the right hand side of (\ref{e.2tensor}) is $$ ({\bf 1}_{x'} \wedge h) \otimes
 f^3 \cdot \mathbf{F}_{x}(v_1, v_2, v_3) \otimes {\bf 1}_{x'} $$ by Lemma \ref{l.parallelIII}. Using the compatibility of the tensor multiplication by $\sA_x^{x'}$ in the identification (iii) of Proposition \ref{p.sD},    this is equal to   $$f^2 \cdot (\mathbf{1}_{x'} \wedge h) \otimes \mathbf{F}_x(v_1, v_2, v_3) \otimes (f \cdot \mathbf{1}_{x'}) $$ $$ = (f \wedge \mathbf{1}_x) \otimes \mathbf{F}_x(v_1, v_2, v_3) \otimes f.$$ Thus (\ref{e.2tensor}) gives \begin{eqnarray*} [f \otimes v_3, [f \otimes v_2, {\bf 1}_x \otimes v_1] ]= - ({\bf 1}_{x_0} \wedge f) \otimes \mathbf{F}_{x}(v_1, v_2, v_3) \otimes f.\end{eqnarray*} This proves (\ref{e.fg2}) when $u_3  $ is a nonconstant function  $f \in U_{x}$.
\end{proof}

        \section{From distributions with symbols $\fg_+(F)$ to nondegenerate lines on contact manifolds}\label{s.zelenko}

\begin{definition}\label{d.Levi}
Let $D \subset TM$ be a distribution on a complex manifold $M$ and let $D^{\perp} \subset T^{\vee} M$ be its annihilator.
\begin{itemize} \item[(i)] Consider the null space of the Levi tensor at $y \in M$,
$${\rm Null}({\rm Levi}_y^D):= \{ v \in D_y \mid {\rm Levi}^D_y (v, u) =0 \mbox{ for all } u \in D_y \}.$$ As $y$ varies over $M$, the null spaces define a distribution ${\rm Ch}(D) $ on a Zariski-open subset of $M$, called the {\em Cauchy characteristic} of $D$. It is easy to see that this is an integrable distribution, defining a holomorphic foliation on the Zariski-open subset of $M$.
\item[(ii)] For each $y \in M$ and  $\alpha \in D_y^{\perp} \subset T^{\vee}_y M,$ define the anti-symmetric form  on $D_y$, $$  \alpha \circ {\rm Levi}_y^D: \wedge^2 D_y \to \C$$ and its null-space  by $${\rm Null}(\alpha \circ {\rm Levi}_y^D) = \{ v \in D_y \mid \alpha \circ {\rm Levi}_y^D (v, u) = 0 \mbox{ for all } u \in D_y\}.$$  
\end{itemize} \end{definition}

\begin{notation}\label{n.vartheta}
The projectivization $\BP T^{\vee} M$ of a complex manifold $M$ has a natural $\sO_{\BP T^{\vee}M} (1)$-valued 1-form $\vartheta^M$ such that ${\rm Ker}(\vartheta^M)$ is a contact structure on $\BP T^{\vee}M$ (see Notation 2.4 of \cite{Hw22} for more details). \end{notation}

The following
 is  Proposition 3.6 of \cite{Hw22} (note that ${\rm Levi}^{\widetilde{H}}_{\alpha}$ is written as $\check{\rm d}_{\alpha} (\vartheta^M|_{\BP D^{\perp}})$ in \cite{Hw22}).

\begin{proposition}\label{p.3.6} In the setting of Definition \ref{d.Levi},
the restriction $\vartheta^M|_{\BP D^{\perp}}$ defines a distribution $\widetilde{\sH}$ of corank 1 on the manifold $\BP D^{\perp}$.
 For any $y \in M$ and any nonzero $a \in D^{\perp}_y$ with the corresponding point $\alpha \in \BP D^{\perp}_y$, we have an isomorphism
$$ {\rm Null} ({\rm Levi}^{\widetilde{\sH}}_{\alpha}) \cong {\rm Null}(\alpha \circ {\rm Levi}_y^D) \subset D_y$$ induced by the natural projection $\BP D^{\perp} \to M$. \end{proposition}

The following is the precise formulation of going from the left hand side to the right hand side of (\textsf{Corr.$m$}) in Section \ref{s.1}.

\begin{theorem}\label{t.lines}
Let $\sD \subset TY$ be a distribution regular at every point of $ Y$ such that the symbol algebra ${\rm symb}_y(\sD)$ is isomorphic to $\fg_+(F_y)$ in Definition \ref{d.fg+} for a nondegenerate cubic form $F_y$ for each $y \in Y$. Let $\sW \subset T^{\vee} Y$ be the annihilator of $\partial \sD$. Denote by $\varrho: \BP \sW \to Y$ the $\BP^1$-bundle obtained by the projectivization of $\sW$. Then for each $y \in Y$, there exists an open neighborhood $O^y \subset Y$ of $y$ and a submersion  $\mu_y: \varrho^{-1}(O^y) \to X^y$ of relative dimension $m$ onto a complex manifold $X^y$ equipped with a contact structure $\sH^y \subset TX^y$  such that
 \begin{itemize}
 \item[(i)] the distribution of corank 1 on $\varrho^{-1} (O^y)$ given by the restriction of  $\vartheta^Y$  to the submanifold $\BP \sW \subset \BP T^{\vee}$ has the fiber at $\alpha \in \varrho^{-1} (O^y) \subset \BP \sW$ equal to   $({\rm d}_{\alpha} \mu_y)^{-1} (\sH^y_{\mu_y(\alpha)})$;
 \item[(ii)]
    for  two distinct points $\alpha_1 \neq \alpha_2 \in \BP \sW_y$, the two vector spaces ${\rm Ker}({\rm d}_{\alpha_1} \mu_y)  \subset T_{\alpha_1}  \BP \sW$ and ${\rm Ker}({\rm d}_{\alpha_2} \mu_y)  \subset T_{\alpha_2}  \BP \sW$ satisfy $${\rm d}_{\alpha_1} \varrho ({\rm Ker}({\rm d}_{\alpha_1} \mu_y))  \ \cap \ {\rm d}_{\alpha_2} \varrho ({\rm Ker}({\rm d}_{\alpha_2} \mu_y)) \ = \ 0;$$  \item[(iii)] the submersion $\mu_y$ sends the fibers of $\varrho$  to   lines on $(X^y, \sH^y)$;
    \item[(iv)] the  lines  on $(X^y, \sH^y)$ in (iii) form a $(3m+2)$-dimensional family, i.e., all of their small deformations on $X^y$   come from the $\mu_y$-images of the fibers of $\varrho$; and \item[(v)] the lines in (iii) are nondegenerate.  \end{itemize}
 \end{theorem}

We need the following two lemmata. We skip the proof of the first one, which is straightforward. The second is well-known (e.g. Lemma 3.5 of \cite{Ke}), but we reproduce the proof for the reader's convenience.

    \begin{lemma}\label{l.varsigma}
In Definition \ref{d.fg+},     denote by $\varsigma: \wedge^2 (\fg_1 + \fg_2) \to \fg_3$ be the homomorphism defined by the Lie bracket of $\fg$ modulo $\fg_1 + \fg_2$. It is independent of $F$ and for any $a \in \fg_3^{\vee}$ the anti-symmetric form $a \circ \varsigma$ on $\fg_1 + \fg_2$ satisfies $${\rm Null}(a \circ \varsigma) = a^{\perp} \otimes V \  \subset \ \fg_1,$$
where $a^{\perp} = \{ u \in \U \mid a (u) =0\}.$
        \end{lemma}

\begin{lemma}\label{l.unbendable}
For a contact manifold $(X, \sH)$ of dimension $2m+3$ with the contact line bundle $\sL = TX/\sH$, let $C \subset X$ be a smooth rational curve such that  $\sL|_C\cong \sO(1)$ and  $TX|_C$ is semipositive, i.e., $$TX|_C \cong \sO(l_1) \oplus \cdots \oplus \sO(l_{2m+3})$$ for some nonnegative integers $ l_1 \geq \cdots \geq l_{2m+3} \geq 0.$
 Then $C$ is a line on $(X, \sH)$. \end{lemma}

\begin{proof}
It suffices to show that $C$ is unbendable.
Let $\lambda: TX \to \sL$ be the quotient homomorphism.
Since $TC \cong \sO(2)$ and $\sL|_C \cong \sO(1)$, we have $\lambda (T C) =0, $
 which implies that $C$ is tangent to $\sH = {\rm Ker}(\lambda)$.
Since $\det TX = L^{\otimes (m+2)}$, we have \begin{equation}\label{e.sum} l_1 + \cdots + l_{2m+3} =  m+2.\end{equation}
Let $$f^*\sH \cong \sO(k_1) \oplus \cdots \oplus \sO(k_{2m+2})$$ for $k_1 \geq \cdots \geq k_{2m+2}$
Since ${\rm Levi}^{\sH}$ gives  symplectic forms on  fibers of $\sH$ with values in $\sL$, we have $$k_1 + k_{2m+2} = k_2 + k_{2m+1} = \cdots = k_m + k_{m+3} = k_{m+1} + k_{m+2} = 1.$$  Thus $k_1 \geq \cdots \geq k_{m+1} \geq 1$ and $0 \geq k_{m+2} \geq \cdots \geq k_{2m+2}$. Combining $2 \leq k_1 \leq l_1$ and $k_i \leq l_i$ for $1 \leq i \leq 2m+2$ with (\ref{e.sum}), we conclude $$l_1 =k_1= 2, \ l_2 = \cdots = l_{m+1} = k_2 = \cdots = k_{m+1} = 1, $$ $$ l_{m+2} = \cdots = l_{2m+3} = k_{m+1} = \cdots k_{2m+1} =0, \  k_{2m+2} = -1.$$ Thus $C$ is unbendable. \end{proof}

    \begin{proof}[Proof of Theorem \ref{t.lines}]
By Lemma \ref{l.varsigma},  for any $y \in Y$, any $\alpha_1 \neq \alpha_2 \in \BP \sW_y$ and corresponding nonzero vectors $a_1 \neq a_2 \in \sW_y$, we have
\begin{equation}\label{e.1} \dim {\rm Null}(a_i \circ {\rm Levi}^{\partial\sD}_y) = m \end{equation} and \begin{equation}\label{e.2} {\rm Null}(a_1 \circ {\rm Levi}^{\partial \sD}_y) \cap  {\rm Null}(a_2 \circ {\rm Levi}^{\partial\sD}_y) = 0.
\end{equation} Applying Proposition \ref{p.3.6} with $M=Y, D = \partial \sD, D^{\perp} = \sW$,
(\ref{e.1}) shows that  the Cauchy characteristic of $\widetilde{\sH}$ is a distribution defined on the whole $\BP (\partial \sD)^{\perp}$.
For each $y \in Y$, we can choose a neighborhood $O^y$ such that  the space of leaves of the foliation ${\rm Ch}(\widetilde{\sH})$ in  $\varrho^{-1} (O^y)$  becomes a complex manifold $X^y$ (e.g. by Lemma 5.6 of \cite{HL}) and the quotient map becomes a submersion  $\mu_y: \varrho^{-1}(O^y) \to X^y$. Then $\widetilde{\sH}$ descends to a contact structure $\sH^y \subset T X^y$ on $X^y$, satisfying (i). (ii) follows from (\ref{e.2}).

To check (iii), let $C \subset X^y$ be the image of  $\varrho^{-1}(y)$ under $\mu_y$. Then $\sL^y|_C \cong \sO(1)$ for the contact line bundle $\sL^y= T X^y/\sH^y$ by (i) together with the fact that $\vartheta^M|_{\BP W}$ is $\sO_{\BP \sW}(1)$-valued 1-form. Since $\mu_y$ is a submersion, the vector bundle $T X^y|_C$ is semipositive. This shows that $C$ is a line by Lemma \ref{l.unbendable}.

To prove (iv), we need to check that there exists no arc $\Delta \subset Y$ with $y \in \Delta$
such that the surface $\varrho^{-1}(\Delta)$ is sent to a single curve in $X^y$ by $\mu_y$. But if such an arc $\Delta$ exists, for two distinct point $\alpha_1 \neq \alpha_2 \in \BP \sW_y$,  we have
$$T_y \Delta \subset {\rm d}_{\alpha_1} \varrho({\rm Ker}({\rm d}_{\alpha_1} \mu_y) \cap {\rm d}_{\alpha_2} \varrho({\rm Ker})({\rm d}_{\alpha_2} \mu_y),$$ a contradiction to (ii).

By (iv), after shrinking $O^y$ if necessary, we may regard $O^y$ as an open subset of
${\rm Lines}(X^y, \sH^y)$ and identify $\varrho^{-1}(O^y)$ with the open subset $Z:= \rho^{-1}(O^y)$ in the universal family.  By Proposition \ref{p.3.6} and Lemma \ref{l.varsigma}, the distribution $\sD$ is spanned by the images of ${\rm d}\varrho ({\rm Ker}({\rm d}_{\alpha} \mu_y))$ for $\alpha \in \varrho^{-1}(y)$. Thus $\sD|_{O^y}$ agrees with the distribution $\sD$ of Definition \ref{d.tensor} under the identification of $O^y$ as an open subset of ${\rm URC}(X^y)$.
Let $C := \mu_y(\varrho^{-1}(y))$ be the line in $(X^y, \sH^y)$ corresponding to $y$.  By Proposition \ref{p.II}, we see that $\II_{\sC_x,z}$ for a point $z \in \varrho^{-1}(y)$ satisfies (i) and (ii) of Definition \ref{d.LM}. By Proposition \ref{p.III}, we see that $\III_{\sC_x, z}$ satisfies (iii) and (iv) of Definition \ref{d.LM}. It follows that $C$ is a nondegenerate line.
\end{proof}

The arguments in the proof of Theorem \ref{t.lines} show that the constructions in Theorem \ref{t.symbol} and Theorem \ref{t.lines} are the inverse of each other.



\bigskip
Jun-Muk Hwang(jmhwang@ibs.re.kr)

\smallskip

Center for Complex Geometry,
Institute for Basic Science (IBS),
Daejeon 34126, Republic of Korea

\bigskip

Qifeng Li(qifengli@sdu.edu.cn)

\smallskip

School of Mathematics,
Shandong University,
Jinan 250100, China

\end{document}